\documentclass[11pt,english]{smfart}
\usepackage[latin1]{inputenc}
\usepackage[T1]{fontenc}
\usepackage[english,french]{}
\usepackage{hyperref}
\usepackage{amssymb,amsmath,url,xspace,smfthm,euscript,enumerate}
\usepackage{fancyhdr}
\usepackage{layout}
\def\notdiv{\,\nmid\,}
\def\div{\,\vert\,}
\def\too{\relbar\lien\rightarrow}

\let\ds=\displaystyle
\let\st=\scriptstyle

  \def\N{\mathbb{N}}
  \def\R{\mathbb{R}}
  \def\C{\mathbb{C}}
  \def\Q{\mathbb{Q}}
  \def\Z{\mathbb{Z}}
  \def\F{\mathbb{F}}
 \def\zz{\mathbb{Z}}
\def\virg{\raise 2pt \hbox{,}\,\,}
\def\rond{{\scriptstyle\circ}}
\def\ab{{\rm ab}}

\def\Cl{{\mathcal C}\hskip-2pt{\ell}}

\def\plus{\ds\mathop{\raise 2.0pt \hbox{$\bigoplus $}}\limits}
\def\mult{\ds\mathop{\raise 2.0pt \hbox{$\bigotimes$}}\limits}
\def\prd{ \ds\mathop{\raise 2.0pt \hbox{$  \prod   $}}\limits}
\def\Cap{ \ds\mathop{\raise 2.0pt \hbox{$\bigcap   $}}\limits}
\def\Cup{ \ds\mathop{\raise 2.0pt \hbox{$\bigcup   $}}\limits}
\def\sm{  \ds\mathop{\raise 2.0pt \hbox{$  \sum    $}}\limits}

\def\w{{\bf w}}
\def\fin{\vbox{\hrule\hbox to 7.2pt{\vrule height 7pt\hfil\vrule}\hrule}}
\def \tensorZp{\otimes{\raise -0.8pt \hbox{\!\!$_{_{\zz_{\!p}}}$}}}
\def \tensorZ{\otimes{\raise -0.8pt \hbox{\!\!$_{_{\zz}}$}}}
\let\st=\scriptstyle  
\let\cal=\mathcal

\def\fin{\vbox{\hrule\hbox to 7.2pt{\vrule height 7pt\hfil\vrule}\hrule}}
\def\lien{\mathrel{\mkern-4mu}}
\def\wt{\widetilde}
\def\ov{\overline} 
\def\wh{\widehat} 

\newenvironment{defis}{\begin{enonce}{Definitions}}{\end{enonce}}
\newenvironment{remas}{\begin{enonce}{Remarks}}{\end{enonce}}

\theoremstyle{plain}
\newcounter{note}

\begin{document}

 \title{On the structure of the  Galois group \\ of  the Abelian closure of a number field}  

\author[Georges  Gras ]{Georges  {\sc Gras}}
\curraddr{Villa la Gardette, chemin Ch\^ateau Gagni\`ere, F-38520 Le Bourg d'Oisans}
\email{g.mn.gras@wanadoo.fr}
\urladdr{http://monsite.orange.fr/maths.g.mn.gras/}

\begin{abstract} Following a paper by Athanasios Angelakis and Peter Stevenhagen on the determination of
imaginary quadratic fields having the same  absolute Abelian Galois group $A$, we study this property for arbitrary number fields. We show that such a property is probably not easily generalizable, apart from imaginary quadratic fields, because of some $p$-adic obstructions coming from the global units. By restriction to the $p$-Sylow subgroups of~$A$, we show that the corresponding study is related to a generalization of the classical notion of $p$-rational fields. However, we obtain some non-trivial information about the structure of the profinite group $A$, for every number field, by application of results published in our book on class field theory.
\end{abstract}

\begin{altabstract} A partir d'un article de Athanasios Angelakis et Peter Stevenhagen sur la d\'etermination de corps quadratiques imaginaires ayant le m\^eme  groupe de Galois Ab\'elien absolu $A$, nous \'etudions cette propri\'et\'e  pour les corps de nombres quelconques. Nous montrons qu'une telle propri\'et\'e n'est probablement pas facilement g\'en\'eralisable, en dehors des corps quadratiques imaginaires, en raison d'obstructions $p$-adiques  provenant des unit\'es globales. En se restreignant aux $p$-sous-groupes de Sylow de $A$, nous montrons que l'\'etude correspondante est li\'ee \`a une g\'en\'eralisation de la notion classique de corps $p$-rationnels. Cependant, nous obtenons des informations non triviales sur la structure du groupe profini $A$, pour tout corps de nombres, par application de r\'esultats publi\'es dans notre livre sur la th\'eorie du corps de classes.
\end{altabstract}

\date{February 18,  2013}

\keywords{Class field theory; Abelian closures of number fields; Abelian profinite groups; $p$-ramification; $p$-rational fields; Group extensions}

\subjclass{Primary 11R37; 11R29; Secondary 20K35}

\maketitle

\vspace{-0.8cm}
\section{Introduction -- Notations} Let $K$ be a number field of signature $(r_1,\,r_2)$ for which
$r_1+2\,r_2 = [K:\Q]$,  and let $A_K$ be the Galois group 
${\rm Gal}(\ov K^{\rm ab}/ K)$ where $\ov K^{\rm ab}$ is the maximal Abelian pro-extension of $K$.
The question that was asked was the following: in what circumstances the groups $A_{K_1}$ and $A_{K_2}$ are isomorphic groups when $K_1$ and $K_2$ are two non-conjugate number fields ? 

\smallskip
A first paper on this subject was published in [O] by M. Onabe. In [AS], A. Angelakis and P. Stevenhagen show that $A_K \simeq \wh \Z^2 \times \prd_{n \geq 1} \Z/n\,\Z$ for a specific family of imaginary quadratic fields.
In this paper, we prove (under the Leopoldt conjecture) that, for any number field $K$, the group $A_K$ {\it contains} a subgroup isomorphic to $\wh \Z^{r_2+1} \! \times \prd_{n \geq 1} \!\! \Big( (\Z/2\,\Z)^\delta  \times\Z/\w \,n \,\Z \Big)$, where $\delta = 1$ if $K \cap \Q(\mu_{2^\infty})$ is a non-trivial extension of $\Q$
distinct from $\Q(\mu_4)$, $\delta = 0$ otherwise, and where $\w = \prd_p \w_p$ is an integer whose local factors $\w_p$ depend simply on the intersections $K\cap\Q(\mu_{p^\infty})$; then we give a class field theory interpretation of the quotient of $A_K$ by this~subgroup, quotient which measures the defect of $p$-rationality for all $p$. 

\smallskip
Such isomorphisms are only isomorphisms of Abelian profinite groups for which
Galois theory  and, a fortiori, arithmetical objects (decomposition and inertia groups)
are not effective.

\smallskip
When an isomorphism is canonical (for instance if it is induced by the reciprocity map of class field theory), we shall write $\, \ds\mathop{\simeq}^{\rm can} \,$ contrary to the non-canonical case denoted
$\ \ds\mathop{\simeq}^{\rm nc}\ $ if necessary.

\smallskip
Let $p$ be a fixed prime number and let 
$$\hbox{  $H$, $H_p$, $H_{\rm ta}$, $\wt K_p$, }$$
be the $p$-Hilbert class field (in ordinary sense), the maximal $p$-ramified (i.e., unramified outside~$p$) Abelian pro-$p$-extension of $K$ (in ordinary sense), the maximal tamely ramified Abelian pro-$p$-extension of $K$ (in restricted sense), the compositum of the $\Z_p$-extensions of $K$, respectively. Then let 
$$\hbox{ ${\mathcal T}_p := {\rm Gal}(H_p/\wt K_p)\ $ and $\ \Cl_p :=  {\rm Gal}(H/ K)$  }$$
canonically isomorphic to the $p$-class group of~$K$.
The groups ${\mathcal T}_p$ and $\Cl_p$ are finite groups.

\smallskip
For any finite place $v$ of $K$, we denote by $K_v$ the corresponding completion\,\footnote{As in [Gr1, I, \S\,2], we consider that $K_v = i_v(K)\,\Q_\ell \subset \C_\ell$ for a suitable embedding $i_v$ of the number field $K$
in $\C_\ell$, where $\ell$ is the residue characteristic.} of $K$, then  by 
$$\hbox{ $U_v := \{u \in K_v, \vert\, u \,\vert_v =1 \}\ $ and $\ U_v^1 := \{u \in U_v, \vert\,u - 1\,\vert_v < 1 \}$, } $$
the unit group and principal unit group of $K_v$, respectively. 
 So, $U_v/U_v^1$ is isomorphic to the multiplicative group of the residue field $F_v$ of $K$ at $v$. We shall denote by   $\ell$ the characteristic of~$F_v$; then $U_v^1$ is a $\Z_\ell$-module.
If $v$ is a real infinite place, we put $K_v = \R$, $U_v = \R^{\times}$, $U_v^1 = \R^{\times +}$, hence $F_v^\times = \{\pm 1\}$, according to [Gr1, I.3.1.2, II.7.1.3].

\smallskip
The structure of $ {\rm Gal}(H_p/K)$ can be summarized by the following diagram, from [Gr1, III.2.6.1, Fig. 2.2]
under the Leopoldt conjecture for $p$ in $K$:
\unitlength=0.6cm
$$\vbox{\hbox{\hspace{-2cm}  \begin{picture}(11.5,5.7)
%
% horizontales
\put(8.4,4.50){\line(1,0){2.5}}
\put(4.3,4.50){\line(1,0){2.5}}
\put(4.8,2.50){\line(1,0){2.1}}

\bezier{350}(4,4.8)(7.5,5.6)(11,4.8)
\put(7.2,5.5){$\st {\mathcal T}_p$}

% vertical
\put(3.50,2.9){\line(0,1){1.20}}
\put(3.50,0.9){\line(0,1){1.20}}
\put(7.50,2.9){\line(0,1){1.20}}
\bezier{350}(3.0,0.8)(2.1,2.5)(3.0,4.2)
\put(1.2,2.4){$\st \Z_p^{r_2+1} $}
% obliques

\bezier{200}(3.9,0.5)(5.7,0.5)(7.2,2.2)
\put(6.4,0.8){$\st \Cl_p$}

\bezier{250}(8.1,2.5)(10.6,2.5)(11.1,4.2)
\put(10.5,2.6){$\st \big(\prod_{v \div p} U_v^1\big) \big/E \otimes\Z_p$}

\put(11.,4.4){$H_p$}
\put(6.9,4.4){$\wt K_{p} H $}
\put(3.2,4.4){$\wt K_{p}$}
\put(7.2,2.4){$H $}
\put(2.8,2.4){$\wt K_{p} \!\cap\! H $}
\put(3.3,0.40){$K$}
\end{picture}   }} $$
\unitlength=1.0cm

where $E$ is the group of global units of $K$ and where $E \otimes\Z_p$ is diagonally embedded
with the obvious map $i_p :=  (i_v)_{v\div p}$.

\smallskip
To characterize the notion of $p$-rationality (see Definition 2.1 and Remarks 2.2),
we shall make use of some $p$-adic logarithms as follows:

\smallskip
(i) We  consider the $p$-adic logarithm ${\rm log}_p  :  K^\times  \too\ds \prd_{v \vert p} K_v$ defined by
${\rm log}^{}_{p} = {\rm log}\,\rond\, i^{}_{p}$ on $K^\times$, where ${\rm log} \ : \ \C_p^\times \too \C_p$ 
is the Iwasawa extension of the usual $p$-adic logarithm.

\smallskip
(ii) We then define the quotient $\Q_p$-vector space
$ {\mathcal L}_{p} := \Big( \prd_{v \vert p} K_v \Big)\Big / \Q_p {\rm
log}^{}_{p} (E)$.
 We have, under the Leopoldt conjecture for $p$ in $K$,
${\rm dim}_{\,\Q_p}^{}({\mathcal L}_{p}) = r_2 + 1$.

\smallskip
(iii) Finally, we denote by $ {\rm Log}_{p}$ the map, from the group
$I_p$ of ideals of $K$ prime to $p$, to ${\mathcal L}_{p}$, sending ${\mathfrak a} \in I_p$ 
to ${\rm Log}_p ({\mathfrak  a})$ defined as follows.
If $m$ is such that ${\mathfrak  a}^m = (\alpha)$ with $\alpha \in K^\times$, we set
${\rm Log}_{p} ({\mathfrak  a}) := \hbox{$\frac{1}{m}$} {\rm log}^{}_{p} (\alpha)
+ \Q_p {\rm log}^{}_{p} (E)$; this does not depend on the choices of $m$ and $\alpha$.

\section{Structure of the Galois group of the Abelian closure of  a number field}

\subsection{Class field theory  -- Fundamental diagram -- $p$-rationality}
Let $\ov K^\ab{\st (p)}$ be the maximal Abelian pro-$p$-extension of~$K$.
In [Gr1, III.4.4.1], we have given  (assuming the Leopoldt conjecture for $p$ in $K$) the following  
diagram for the structure of ${\rm Gal}(\ov K^\ab{\st (p)} / K)$, isomorphic to the $p$-Sylow subgroup of $A_K$:
\unitlength=0.65cm
$$\vbox{\hbox{\hspace{-1.5cm}  \begin{picture}(11.5,6.5)
%
% horizontales
\put(8.5,4.50){\line(1,0){3.0}}
\put(1.5,4.50){\line(1,0){5.9}}
\put(1.5,2.50){\line(1,0){5.9}}
%
% verticales
\put(1.0,2.9){\line(0,1){1.20}}
\put(1.0,1.4){\line(0,1){0.70}}

\put(8.00,2.9){\line(0,1){1.20}}
\bezier{400}(1.2,4.9)(6.45,6.1)(11.7,4.9)
\put(5.6,5.8){$\st {\prod_{v \notdiv p}{ ( F_v^\times)_p  }}$}

\bezier{280}(8.8,2.5)(11.5,2.8)(11.9,4.2)
\put(11.,2.5){$\st \prod_{{v\,|\,p}}{U_v^{1}} $}

\put(9.2,4){$\st  E \otimes\zz_p$}

\put(11.7,4.4){$\ov K^\ab{\st (p)}$}

\put(7.65,4.4){$M_p$}
\put(0.7,4.4){$H_p$}

\put(7.55,2.4){$H_{\rm ta}$}
\put(0.7,2.4){$H$}

\put(0.7,0.90){$K$}

\end{picture}   }} $$
\unitlength=1.0cm

\vspace{-0.4cm}
where $( F_v^\times)_p$ is the $p$-Sylow subgroup of the  multiplicative group of the residue field 
$F_v$ of $K$ at $v$. This also concerns the real infinite places for which $F_v^\times = \{ \pm 1 \}$.
 In this diagram, $M_p$ is the direct compositum, over $H$, of $H_p$ and $H_{\rm ta}$.

\smallskip
The diagonal embeddings $i_{\rm ta} := (i_v)_{v\notdiv p}$ and $i_p := (i_v)_{v\div p}$ of $E \otimes\zz_p$ in ${\prd_{v \notdiv p}{ ( F_v^\times)_p}}$
and $\prd_{{v\,|\,p}}{U_v^{1}}$, respectively, are injective (under the Leopoldt conjecture for the second one).

\smallskip
 Let ${\mathcal U}_p := 
{\prd_{v \notdiv p}{ ( F_v^\times)_p}} \times \prd_{{v\,|\,p}}{U_v^{1}}$ be the $p$-Sylow subgroup of the group of unit id\`eles ${\mathcal U} := \prd_v U_v$ of $K$, and let  $\rho$ be the reciprocity map on~${\mathcal U}_p$. 

\smallskip
The kernel of $\rho$  is 
$i(E \otimes\zz_p)$, where $i=(i_{\rm ta}, i_p)$; this yields ${\rm Gal}(\ov K^{\rm ab}{\st (p)}/H)
\simeq {\mathcal U}_p/ i(E \otimes\zz_p)$, ${\rm Gal}(\ov K^{\rm ab}{\st (p)}/H_p) = 
\rho \Big ( {\prd_{v \notdiv p}{ ( F_v^\times)_p}} \times \{ 1\} \Big ) \simeq {\prd_{v \notdiv p}{ ( F_v^\times)_p}}$ since $\Big ({\prd_{v \notdiv p}{ ( F_v^\times)_p}} \times \{ 1\} \Big)\ \cap\  i\,(E \otimes\zz_p) = 1$, and similarly 
${\rm Gal}(\ov K^{\rm ab}{\st (p)}/H_{\rm ta}) = \rho\Big ( \{ 1\} \times \prd_{{v\,|\,p}}{U_v^{1}} \Big ) \simeq \prd_{{v\,|\,p}}{U_v^{1}}$ (see [Gr1, III.4.4.5.1]). 

This will be useful in Section 3.

\begin {defi}  The number field $K$ is said to be $p$-rational (see  [MN], [GJ], [JN], and [Gr1, IV, \S\,b, 3.4.4 ]) if it satisfies the Leopoldt conjecture for $p$ and if ${\mathcal T}_p = 1$.
\end{defi}

\begin{remas} {\rm Assuming the Leopoldt conjecture for $p$ in $K$, we have:
 
\smallskip
(i) From [Gr1, IV.3.4.5], the $p$-rationality of $K$ is equivalent to the following three conditions:

\smallskip
$\ \ \bullet\ $ $\prd_{v|p} \mu_p(K_v) = i_p(\mu_p(K))$, where $\mu_p(k)$ denotes, for any field $k$, the group of roots of unity of $k$ of $p$-power order, 

\smallskip
$\ \ \bullet\ \ $the $p$-Hilbert class field $H$ is contained in the compositum $\wt K_p$ of the 
$\Z_p$-extensions of~$K$; this is equivalent to $\Cl_p \, {\ds\mathop{\simeq}^{\rm can} }\,
 \Z_p {\rm Log}_p(I_p)\Big / \Big(\prd_{v|p} {\rm log}(U_v^{1}) + \Q_p {\rm log}_p(E) \Big)$, 
which can be non-trivial,

$\ \ \bullet\ $ $\Z_p{\rm log}_p(E) \hbox{ is a direct summand in the $\Z_p$-module } \prd_{v|p}
{\rm log}(U_v^{1})$, which expresses the minimality of the valuation of the $p$-adic regulator.

\smallskip
(ii) For a $p$-rational field $K$, we have 
${\rm Gal}(\ov K^\ab{\st (p)}/K) \, {\ds\mathop{\simeq}^{\rm nc}} \,  \Z_p^{r_2 +1} \times \prd_{v \notdiv p} (F_v^\times)_p$, with (canonically) ${\rm Gal}(\ov K^\ab{\st (p)}/\wt K_p)\, {\ds\mathop{\simeq}^{\rm can} }\,
\prd_{v \notdiv p} ( F_v^\times)_p$.

\smallskip
(iii) Let $\wt {\mathcal K}_\infty$ be the compositum of the $\wt K_p$. By assumption, it is
the maximal $\wh\Z$-extension of~$K$ for which 
${\rm Gal}(\wt {\mathcal K}_\infty/K) \, {\ds\mathop{\simeq}^{\rm nc} }\,  \wh\Z^{r_2 +1} $. 
A sufficient condition to get ${\rm Gal}(\ov K^\ab/K) \, {\ds\mathop{\simeq}^{\rm nc} }\, 
\wh\Z^{r_2 +1} \times \prd_{v } F_v^\times$ is that $K$ be $p$-rational for all~$p$. 
The stronger condition $H_p = \wt K_p$ for all $p$ (i.e., $p$-rationality of $K$ for all $p$) is
equivalent to the class field theory isomorphism 
${\rm Gal}(\ov K^\ab/\wt {\mathcal K}_\infty) \, {\ds\mathop{\simeq}^{\rm can} }\,\prd_{v } F_v^\times$. }
\end{remas}

\subsection{Structure of $\prd_v F_v^\times$}
Let $(F_v)_v$ be the family of residue fields of $K$ for its finite or real infinite places $v$.
We intend to give, for all prime numbers $p$, the structure of the $p$-Sylow subgroup of $\prd_v F_v^\times$.
If $v\div p$, then $(F_v^\times)_p =1$; so we can restrict ourselves
to $\prd_{v\notdiv p} (F_v^\times)_p$. 

We shall prove that there exist  integers $\delta \in \{0, 1\}$ and $\w \geq 1$ such that
$$\prd_{v\notdiv p} ( F_v^\times)_p \,\ds\mathop{\simeq}^{\rm nc} \, \Big( \prd_{n \geq 1}
\Big( (\Z/2\,\Z)^\delta  \times \Z/\w \, n \,\Z\Big) \Big)_p . $$

 The property giving such an isomorphism is that for any given $p$-power $p^k$, $k\geq 1$, the two $p$-Sylow subgroups have $0$ or infinitely many (countable) cyclic direct components of order~$p^k$. 

\smallskip
It is obvious that for any $p$, $\Big( \prd_{n \geq 1} \Big( (\Z/2\,\Z)^\delta  \times \Z/\w \, n \,\Z\Big) \Big)_p$ has  $0$ or infinitely many cyclic direct components of order $p^k$ for any $k\geq 1$; 
moreover, in $\Big( \prd_{n \geq 1} \big( \Z/\w\, n \,\Z\big) \Big)_p$ there is no  direct component of order $p^k$, $k\geq 1$, if and only if $p^{k+1} \div \w $.

\begin{rema} {\rm Write $\w  = \prd_{q\  {\rm prime}} q^{\lambda_q}$ and put $\w^1 := \prd_{\lambda_q\geq 2} q^{\lambda_q}$ so that $\w^0 := \w /\w^1$ and $\w^1$ are coprime integers. Then in the above expressions
we can replace $\w$ by $\w^1$.
Indeed, in the two writings $\prd_{n \geq 1} \Z/ \w^0 .\w^1 n \,\Z$ and 
$\prd_{n \geq 1} \Z/\w^1 n \,\Z$, for all $q \div \w^0$ the direct summands of order $q$ are infinite in number.
Then we can ensure that $\w $ will be defined in such a way that $\w^0=1$. }
\end{rema}

\begin{defis} {\rm
(i)  We denote by $\mu(K)$ (resp. $\mu_p(K)$) the group of roots of unity of $K$ (resp. of $p$-power order) 
and for any $e\geq 1$ we denote by $\mu_e$ the group (of order~$e$) of $e$th roots of unity in a field of
characteristic $0$ or $\ell$ prime to $e$.

\smallskip
(ii) Let $Q_\nu$, $\nu\geq 1$, be for any $p$ the unique subfield, of degree $p^\nu$ over $\Q$, 
of the cyclotomic $\Z_p$-extension of~$\Q$.
Then let $Q'_\nu$, $\nu\geq 1$, be for $p=2$ the unique non-real subfield of $\Q(\mu_{2^\infty})$ 
of degree $2^\nu$ over $\Q$.  We put $Q'_0 = Q_0 = \Q$ in any case.}
\end{defis}

\subsubsection{Analysis of the case $p\ne 2$} 
In the study of $\prod_{v\notdiv p} (F_v^\times)_p$,
we can assume that $\vert F_v^\times \vert \equiv 0 \pmod p$, which is equivalent to the
splitting of $v$ in $K(\mu_p)/K$ (this includes the case where $K$ contains $\mu_p$).

\smallskip
{\bf a)} If $K$ contains $\mu_p$, then $\mu_p(K) = \mu_{p^{\nu+1}}$ where $\nu \geq 0$ is the maximal 
integer such that $Q_\nu \subseteq K$, and we get necessarily $\vert F_v^\times \vert \equiv 0 \pmod {p^{\nu+1}}$
for all these places.
We obtain the following tower of extensions (where $\subset$ means a stric inclusion)
$$K =  K(\mu_{p^{\nu+1} } ) \subset  K(\mu_{p^{\nu+2} } )  \subset \cdots $$

From Chebotarev's theorem, for any $m\geq \nu+1$ there exist infinitely many places $v$ of $K$, 
whose inertia group in $K(\mu_{p^{m+1} } ) /K$ is
${\rm Gal}(K(\mu_{p^{m+1} } )/K(\mu_{p^{m} } ))$, cyclic of order $p$; so we get $\vert F_v^\times \vert \equiv 0 \pmod {p^m}$ and $\vert F_v^\times \vert \not\equiv 0 \pmod {p^{m+1}}$ for these places.

\smallskip
{\bf b)} If $K$ does not contain $\mu_p$ and if $K \cap \Q(\mu_{p^\infty}) = Q_\nu$, $\nu \geq 0$,
 we have the tower of extensions
$$K \subset K(\mu_p) = \cdots =  K(\mu_{p^{\nu+1} } ) \subset  K(\mu_{p^{\nu+2} } )  \subset \cdots $$

Since by assumption the places $v$ considered are split in $K(\mu_p)/K$, we still have $\vert F_v^\times \vert \equiv 0 \pmod {p^{\nu+1}}$.
From Chebotarev's theorem,  for any $m\geq \nu+1$ there exist infinitely many places $v$, whose inertia group in $K(\mu_{p^{m+1} } )/K$ is 
${\rm Gal}(K(\mu_{p^{m+1} } )/K(\mu_{p^{m} } ))$, cyclic of order $p$; thus, $\vert F_v^\times \vert \equiv 0 \pmod {p^m}$ and $\vert F_v^\times \vert \not\equiv 0 \pmod {p^{m+1}}$ for these places.

\smallskip
In conclusion, the case $p\ne 2$ leads to identical results from the knowledge of the integer $\nu$.

\subsubsection{Analysis of the case $p= 2$}
 In that case, we always have $\vert F_v^\times \vert \equiv 0 \pmod 2$  in the study of $\prod_{v\notdiv 2} (F_v^\times)_2$ ($v$ finite or real infinite). So we consider $K(\mu_4)/K$.

\smallskip
{\bf a)}  If  $K$ contains $\mu_4$ and if $K \cap  \Q(\mu_{2^\infty}) =: \Q (\mu_{4.2^{\nu} } )$, $\nu \geq 0$,
we have $\vert F_v^\times \vert \equiv 0 \pmod {4.2^\nu}$ for all places,
and the tower of extensions
$$K = K(\mu_{4.2^{\nu} } ) \subset K(\mu_{4.2^{\nu+1} } ) \subset \cdots $$

From Chebotarev's theorem,  for any $m\geq \nu$ there exist infinitely many places $v$ 
whose inertia group, in $K(\mu_{4.2^{m+1} } )/K$,
 is ${\rm Gal}(K(\mu_{4.2^{m+1} } )/K(\mu_{4.2^{m} } ))$, cyclic of order $2$; so 
$\vert F_v^\times \vert \equiv 0 \pmod {4.2^m}$ and
$\vert F_v^\times \vert \not\equiv 0 \pmod {4.2^{m+1}}$ for these places.

\smallskip
{\bf b)}  If  $K$ does not contain $\mu_4$,  we have two possible towers
depending on $K \cap  \Q(\mu_{2^\infty})$:
\begin{eqnarray*}
&&\bullet\  K \cap \Q(\mu_{2^\infty}) =  \Q   : \hspace{3.cm}  K  \subset K(\mu_4) \subset  K(\mu_8 ) \subset \cdots \\
&&\bullet\   K \cap \Q(\mu_{2^\infty}) \in\{ Q_{\nu},  Q'_{\nu}\} ,\ \nu \geq 1  : \hspace{0.6cm}
 K  \subset K(\mu_4) = \cdots =  K(\mu_{4.2^{\nu} } ) \subset  K(\mu_{4.2^{\nu+1} } )  
\subset \cdots 
\end{eqnarray*}

(i) In the first case ($\nu = 0$),  for any $m\geq 1$ Chebotarev's theorem gives infinitely 
many places $v$ whose inertia group in $K(\mu_{2^{m+1} } )/K$ is
${\rm Gal}(K(\mu_{2^{m+1} } )/K(\mu_{2^{m} } ))$, cyclic of order $2$; so we get  $\vert F_v^\times \vert \equiv 0 
\pmod {2^m}$ and $\vert F_v^\times \vert \not\equiv 0 \pmod {2^{m+1}}$ for these places (the real infinite places
are solution, taking $m=1$).

\smallskip
(ii) In the second case ($\nu \geq 1$), we will have two disjoint sets of places of $K$ for the structure of the product 
$\prod_{v \notdiv 2} (F_v^\times)_2$:

\smallskip
-- There exist infinitely many places $v$ inert in $K(\mu_4)/K$. Then $\vert F_v^\times \vert \equiv 0 \pmod {2}$ and $\vert F_v^\times \vert \not\equiv 0 \pmod {4}$ for these places (this includes the real infinite places, if any).

\smallskip
--  For any $m\geq \nu$ ($\nu \geq 1$), Chebotarev's theorem gives infinitely many places $v$ 
whose inertia group in $K(\mu_{4.2^{m+1} } )/K$ is 
${\rm Gal}(K(\mu_{4.2^{m+1} } )/K(\mu_{4.2^{m} } ))$, cyclic of order $2$; a fortiori, these places are split in $K(\mu_4)/K$. So we get $\vert F_v^\times \vert \equiv 0 \pmod {4.2^m}$ and 
$\vert F_v^\times \vert \not\equiv 0 \pmod {4.2^{m+1}}$. 

\smallskip
We see that in the exceptional case $K \cap \Q(\mu_{2^\infty}) \in\{ Q_{\nu},  Q'_{\nu}\}$ with $\nu \geq 1$,
we have a group isomorphism of the form $\prd_{v\notdiv 2} (F_v^\times)_2\simeq \big( \Z/2\,\Z \big)^\N  \times  \prd_{m\geq \nu}\Z/4.2^{m}  \,\Z$. 

\begin{defi} From the above discussion about the number field $K$, we define the integers $\delta \in \{0, 1\}$ and $\w := \prd_{p} \w_p$, where $\w_p$, depending on $K \cap \Q(\mu_{p^\infty})$, is given as follows:

(i) Case $p\ne 2$. Let $\nu\geq 0$ be the maximal integer such that
$Q_\nu \subseteq K$ (thus $\mu_{p^{\nu+1}}$ is the maximal group of roots of unity of $p$-power order 
contained in $K(\mu_p)$ whether  $K$ contains $\mu_p$ or not); we put $\w_p = p^{\nu+1}$ if $\nu \geq 1$ and $\w_p =1$ otherwise (from the use of Remark 2.3).

\smallskip
(ii) If, in the case $p=2$, $K$ contains $\mu_4$, we put $\w_2 = 4.2^\nu$,  where $\nu \geq 0$ is the
 maximal integer such that $Q_\nu \subseteq K$ (in this case, the reasonning with $Q'_\nu$ gives the same integer $\nu$).

\smallskip
(iii)  If, in the case $p=2$,  $K$ does not contain $\mu_4$, let $\nu \geq 0$ be the maximal integer such that
$Q_\nu \subseteq K$ or $Q'_\nu \subseteq K$ (thus $\mu_{4.2^\nu}$ is the maximal
group of roots of unity of $2$-power order of $K(\mu_4)$);  we put $\w_2 =4. 2^{\nu}$
if $\nu \geq 1$ and $\w_2 = 1$ otherwise (from the use of Remark 2.3).

\smallskip
(iv) We put $\delta = 1$ in the case (iii) when $\nu \geq 1$, and $\delta = 0$ otherwise.
\end{defi}

We can state the following result correcting an error discovered by Peter Stevenhagen in the first draft of [AS, Lemma 3.2] as well as in [O] and in the previous versions of our paper reproducing this Lemma;  this will also be corrected in the final paper [AS] in the proceedings volume of ANTS-X, San Diego 2012. 
We refer to  Definitions 2.4 and 2.5 giving $\delta$ and $\w$.

\begin{prop} Let $K$ be a number field.  We have a group isomorphism of the form
$$\prd_v F_v^\times \,\ds\mathop{\simeq}^{\rm nc} \,\prd_{n \geq 1}
\Big( (\Z/2\,\Z)^\delta  \times \Z/\w \, n \,\Z\Big) . $$

We have $\delta = 1$ if and only if $K$ does not contain $\mu_4$ and $8\div \w$.

When $\w = 1$ (the most usual case), then $\delta = 0$ and 
$\prd_v F_v^\times \,\ds\mathop{\simeq}^{\rm nc} \,\prd_{n \geq 1} \Z/n\,\Z$.
\end{prop}

\subsubsection{Examples}
(i) Example with $p=3$. Let $K$ be the maximal real subfield of $\Q(\mu_9)$; 
we have $\w = 9$.
The prime $\ell = 5$ is totally inert in $\Q(\mu_9)$; then for $v\div \ell$, $F_v$ does not contain $\mu_3$
since $\ell^3 = 125 \not \equiv 1 \pmod 3$. But for $\ell = 7$, inert in $K$ and 
split in $\Q(\mu_3)$, we get $F_v^\times = \F_{343}^\times$ which contains $\mu_9$ as expected. 

\smallskip
(ii) Examples with $p=2$. For $K = \Q(\sqrt 2\,)$, we have $\delta=1$ and $\w = 8$. 
The prime $\ell = 7$ splits in $K$ and is inert in $\Q(\mu_4)$; so  for $v\div \ell$, $F_v = \F_7$ 
does not contain $\mu_4$. But for the prime $\ell = 5 \equiv 1 \pmod 4$, inert in $K$
and split in  $\Q(\mu_4)$, we get $F_v= \F_{25}$ which contains $\mu_8$.
For $K = \Q(\sqrt 2\,)$, we get the extra factor $(\Z/2\,\Z)^\N$ and there do not exist any 
cyclic direct component of order $4$.

\smallskip
For $K=\Q(\mu_4)$, we have $\w = 4$ and $F_v = \F_\ell$ ($\ell \equiv 1 \pmod 4$) 
or $F_v = \F_{\ell^2}$  ($\ell \equiv -1 \pmod 4$); so the $2$-Sylow subgroup of $F_v^\times$ 
is at least of order $4$.

\subsection{Consequences for the structure of ${\rm Gal}(\ov K^\ab /K)$}
From Proposition 2.6 and the fundamental diagram (Subsection 2.1), we can state:

\begin{prop} Let ${\mathcal H}_\infty$ be the compositum of the fields $H_p$ (maximal $p$-ramified Abelian pro-$p$-extensions of $K$) for all prime numbers $p$. Then, under the Leopoldt conjecture in $K$ for all~$p$, 
we have a group isomorphism of the form ${\rm Gal}(\ov K^\ab /{\mathcal H}_\infty) \, 
\ds\mathop{\simeq}^{\rm nc} \, \prd_{n \geq 1} \Big( (\Z/2\,\Z)^\delta  \times \Z/\w \, n \,\Z\Big)$.

If $\w = 1$, then  ${\rm Gal}(\ov K^\ab /{\mathcal H}_\infty) \, 
\ds\mathop{\simeq}^{\rm nc} \, \prd_{n \geq 1} \Z/ n \,\Z$.
\end{prop}

We have obtained the following globalized diagram (under the Leopoldt conjecture for all~$p$), 
where ${\mathcal H}_{\rm ta}$ (compositum of the $H_{\rm ta}$) is the maximal Abelian tamely ramified extension 
of~$K$ and ${\mathcal M}_\infty = {\mathcal H}_\infty {\mathcal H}_{\rm ta}$ (direct compositum over 
the Hilbert class field ${\mathcal H}$):
\unitlength=0.65cm
$$\vbox{\hbox{\hspace{-1.5cm}  \begin{picture}(11.5,6.5)
% horizontales
\put(8.5,4.50){\line(1,0){3.0}}
\put(1.8,4.50){\line(1,0){5.4}}
\put(1.8,2.50){\line(1,0){5.4}}
% verticales
\put(1.0,2.9){\line(0,1){1.20}}
\put(1.0,1.4){\line(0,1){0.70}}

\put(8.00,2.9){\line(0,1){1.20}}
\bezier{400}(1.2,4.9)(6.45,6.1)(11.7,4.9)
\put(4.2,6.){$\st {\prod_{n\geq 1 }(  (\Z/2\,\Z)^\delta  \times \Z/\w \, n \,\Z) }$}
\bezier{280}(8.8,2.5)(11.5,2.8)(11.9,4.2)
\put(11.2,2.8){$\st \prod_{v}{U_v^{1}} $}
\put(9.5,4){$\st  E\otimes {\wh\zz}$}
\put(11.8,4.4){$\ov K^\ab$}
\put(7.4,4.4){${\mathcal M}_\infty$}
\put(0.7,4.4){${\mathcal H}_\infty$}
\put(7.7,2.4){${\mathcal H}_{\rm ta}$}
\put(0.7,2.4){${\mathcal H}$}
\put(0.7,0.90){$K$}
\end{picture}   }} $$
\unitlength=1.0cm

\vspace{-0.4cm}
Let ${\mathcal F}_\infty$ be the compo\-situm of some finite extensions $F_p$ of $K$ such that $H_p = \wt K_p\,F_p$
(direct compositum over~$K$). When they are non-trivial, the extensions $F_p/K$ are non-unique $p$-ramified extensions.  We then have 
${\rm Gal}(F_p/K) \, \simeq \, {\mathcal T}_p\ $ and $\ {\rm Gal}({\mathcal F}_\infty/K) \, \simeq\, \prd_p {\mathcal T}_p$. 

 The extension ${\mathcal F}_\infty/K$ is in general non-canonical and conjecturally infinite; its Galois group measures a mysterious degree of complexity of  ${\rm Gal}(\ov K^\ab/K)$; it is trivial if and only if $K$ is $p$-rational for all $p$ (Remark 2.2 (ii)). But we have 
${\rm Gal}({\mathcal H}_\infty/\wt {\mathcal K}_\infty)\  \ds\mathop{\simeq}^{\rm can} \ \prd_p {\mathcal T}_p$
 and ${\rm Gal}({\mathcal H}_\infty/K)\  \ds\mathop{\simeq}^{\rm nc} \  \wh \Z^{r_2+1} 
 \times \prd_p {\mathcal T}_p$.

\vspace{-0.2cm}
\begin{theo} Let $K$ be a number field and let $\ov K^\ab$ be the maximal Abelian pro-extension of~$K$. 
We assume that the $p$-adic Leopoldt conjecture is verified in $K$ for all prime number $p$.

\smallskip
Then there exists an Abelian extension ${\mathcal F}_\infty$ of $K$, with ${\rm Gal}({\mathcal F}_\infty/K)
 \, \ds\mathop{\simeq}^{\rm can} \,  \prd_p {\mathcal T}_p$,
such that ${\mathcal H}_\infty$ is the direct compositum of ${\mathcal F}_\infty$ and the maximal $\wh \Z$-extension~$\wt {\mathcal K}_\infty$ of~$K$, and such that
$${\rm Gal}(\ov K^\ab/{\mathcal F}_\infty) \, \ds\mathop{\simeq}^{\rm nc} \,  \wh \Z^{r_2+1}  \times \prd_{n \geq 1} \Big(  (\Z/2\,\Z)^\delta  \times \Z/\w \, n \,\Z\Big), $$

\vspace{-0.2cm}
  with ${\rm Gal}(\ov K^\ab /{\mathcal H}_\infty) \, \ds\mathop{\simeq}^{\rm nc} \, \prd_{n \geq 1} \Big(  (\Z/2\,\Z)^\delta  \times \Z/\w \, n \,\Z\Big)$, 
where $\delta$, $\w $ are defined in Definition 2.5.

 If $\w = 1$,  we have  a group isomorphism  of the form
${\rm Gal}(\ov K^\ab/{\mathcal F}_\infty) \, \ds\mathop{\simeq}^{\rm nc} \, \wh \Z^{r_2+1} \times \prd_{n \geq 1} \Z/n \,\Z$,  with ${\rm Gal}(\ov K^\ab /{\mathcal H}_\infty) \, \ds\mathop{\simeq}^{\rm nc} \, \prd_{n \geq 1} \Z/ n \,\Z$.
\end{theo}

\vspace{-0.4cm}
\begin{coro} The Galois groups ${\rm Gal}(\ov K^\ab/{\mathcal F}_\infty)$ (up to non-canonical isomorphisms)
are independent of the number fields $K$ as soon as these fields satisfy the Leopoldt conjecture for all $p$,
 have the same number $r_2$ of complex places and the same parameters $\delta, \w$.

\smallskip
Thus, for all totally real number fields $K$ (satisfying the Leopoldt conjecture for all $p$) which do not contain
$\sqrt {2}\,$, we have ${\rm Gal}(\ov K^\ab/{\mathcal F}_\infty) \, \ds\mathop{\simeq}^{\rm nc} \, \wh \Z 
\times \prd_{n \geq 1} \Z/n \,\Z$.
\end{coro}

Of course, the groups ${\rm Gal}({\mathcal F}_\infty/K)$ strongly depend on $K$, even if the parameters $r_2, \delta,\w$ are constant.
 From Remark 2.2 (i), we see that the first two conditions of $p$-rationality involve a finite number of primes $p$, but that the third condition is the most ugly.\,\footnote{For instance, for $K = \Q(\sqrt 2\,)$, the third condition is not satisfied for $p=  13, 31, 1546463, \ldots$ and perhaps for infinitely many primes $p$ depending on Fermat quotients of the fundamental unit [Gr1, III.4.14]. Note that
from [Gr2,~III] or [Gr1, IV.3.5.1], for $p = 2$, the $2$-rational Abelian $2$-extensions of $\Q$ are the subfields of the fields $\Q(\mu^{}_{2^\infty}) \Q(\sqrt{-\ell}\,)$, $-\ell \equiv 5 \pmod 8$, or of the fields
$\Q(\mu^{}_{2^\infty}) \Q\Big (\hbox {$\sqrt{\sqrt{\ell}\ \frac{a-\sqrt{\ell}}{2}}$}\, \Big), \ \ell = a^2+ 4 b^2 \equiv 5 \pmod 8$.} 

\smallskip
So, we are mainly concerned with the imaginary quadratic fields, studied in [AS], for which the third condition is empty; the first and second ones can be verified (for all $p$) probably for infinitely many imaginary quadratic fields as suggested in  [AS, Conjecture 7.1].

\section{A generalization of $p$-rationality}

As we shall see now, we can strengthen a few the previous results about
the first condition involved in the definition of $p$-rationality, condition 
which can be removed for all number~fields.

\smallskip
This concerns the finite $p$-groups $\big(\prod_{v \div p} \mu_p(K_v) \big) \big /i_p(\mu_p(K))$ whose globalization 
 measures the gap between the regular and Hilbert kernels in ${\rm K}_2(K)$  (see [Gr1, II.7.6.1]). 

\smallskip
For all finite place $v$ of $K$ we have $\mu(K_v) \simeq F_v^\times \times \mu^1_v$, where $\mu^1_v$
is the torsion subgroup of $U_v^1$ (it is a finite $\ell$-group where $\ell$ is the residue characteristic);
if $v$ is real infinite, we have $F_v^\times = \{\pm 1\}$ and $\mu^1_v=1$. 

\smallskip
The places (finite in number) such that $\mu^1_v\ne 1$ are called the irregular places of $K$.

\smallskip
We  have $\mu_p(K_v) = \mu^1_v$ if and only if $v \div p$ and $\mu_p(K_v)  \simeq (F_v^\times)_p$ if and only if
$v\notdiv p$.  Let 
$$\Gamma_p := \prd_{v \notdiv p} ( F_v^\times)_p \times \prd_{v \div p} \mu^1_v\simeq \prd_v \mu_p(K_v) . $$

To study the influence of the cyclic factors $\mu^1_v = \mu_p(K_v)$ for $v\div p$, on $\prd_{v \notdiv p} (F_v^\times)_p$, we refer to Definition 2.5 for the definitions of $\nu$, $\delta$, $\w_p$, and to Proposition 2.6.

\smallskip
(i) Case $p\ne 2$.   If $K$  contains $\mu_p$, then $\w_p = p^{\nu+1} = \vert \mu_p(K)\vert$ divides $\vert \mu_p(K_v)\vert$; so the cyclic factor $\mu_p(K_v)$ does not modify the structure.

\smallskip
 If $K$ does not contain $\mu_p$,  we have only to look at the case $\nu\geq 1$ for which
 $\w_p = p^{\nu+1}$. If $\mu_p(K_v)$ is non-trivial ($v \div p$ 
is split in $K(\mu_p)$), $\vert \mu_p(K_v) \vert$ is a multiple of $ p^{\nu+1}$, giving the result.

\smallskip
(ii) Case $p = 2$. If $K$ contains $\mu_4$, $\w_2 = 4.2^\nu = \vert \mu_2(K)\vert$ divides
$\vert\mu_2(K_v)\vert$, hence the result.

\smallskip
  If $K$ does not contain $\mu_4$, we have only to consider the case $K \cap \Q(\mu_{2^\infty})  \in\{ Q_{\nu},  Q'_{\nu}\} ,\ \nu \geq 1$. Then $\delta = 1$ and $\w_2=4.2^\nu$; so $\mu_2(K_v) = \mu_2$ 
(if $v \div 2$ is not split in $K(\mu_4)$) or $\mu_{4.2^m}$, $m\geq \nu$ (if $v$ is split in $K(\mu_4)$), hence the result.

\smallskip
We then have 
$$\Gamma_p \, \ds\mathop{\simeq}^{\rm nc} \, \Big( \prd_{n \geq 1} \Big(  (\Z/2\,\Z)^\delta \times \Z/\w \, n \,\Z\Big) \Big)_p. $$

\smallskip
 Let $H_p^1$ be the subfield of $H_p$ fixed by $\rho(\Gamma_p)$, where $\rho$ is the reciprocity map on 
the $p$-Sylow subgroup ${\mathcal U}_p := \prd_{v \notdiv p} ( F_v^\times)_p \times \prd_{v \div p} U_v^1 \supset \Gamma_p$ of the group of unit id\`eles of $K$. The kernel of $\rho$ is $i (E \otimes\zz_p)$ (see   Subsection 2.1).

\smallskip
We consider $\rho(\Gamma_p) = {\rm Gal}(\ov K^\ab{\st (p)}/H^1_p)$.
Then from the local-global characterization of the Leopoldt conjecture at $p$ (see [Ja, \S\,2.3] or [Gr1, III.3.6.6]),
we get (omitting the embedding~$i$)
$$\rho(\Gamma_p) \, \ds\mathop{\simeq}^{\rm can} \, \Gamma_p/(E \otimes\zz_p) \cap \Gamma_p = 
\Gamma_p/ \mu_p(K). $$

Taking, as in [AS, Lemmas 3.3,  3.4], $v_0$ such that the residue image of $\mu_p(K)$ 
is equal to $( F_{v_0}^\times)_p$, we still get ${\rm Gal}(\ov K^\ab{\st (p)}/H^1_p) \, \ds\mathop{\simeq}^{\rm can} \, \Gamma_p/ \mu_p(K) \, \ds\mathop{\simeq}^{\rm nc} \,  \Big( \prd_{n \geq 1} \Big(  (\Z/2\,\Z)^\delta  \times \Z/\w \, n \,\Z\Big) \Big)_p$.

We note that ${\rm Gal}(H_p/H_p^1) \, {\ds\mathop{\simeq}^{\rm can}} \, \big(\prod_{v\div p } \mu^1_v\big) \big / \mu_p(K)$ (see also [Gr1, III.4.15.3]).

\smallskip
Of course, if ${\mathcal H}^1_\infty \subseteq  {\mathcal H}_\infty$ is the compositum of the $H^1_p$, the globalization gives
$$\hbox{ ${\rm Gal}(\ov K^\ab/{\mathcal H}^1_\infty) \, {\ds\mathop{\simeq}^{\rm nc}} \, \prd_{n \geq 1} \Big(  (\Z/2\,\Z)^\delta  \times \Z/\w \, n \,\Z\Big)$, 
with ${\rm Gal}({\mathcal H}_\infty/{\mathcal H}^1_\infty)\, {\ds\mathop{\simeq}^{\rm can}} \, \Big(\prd_{v} \mu^1_v \Big) \Big / \mu(K)$. } $$

\vspace{-0.1cm}
In other words we have obtained (to be compared with Theorem 2.8 using the extension ${\mathcal F}_\infty$):

\begin{theo} Let $K$ be a number field and let $\ov K^\ab$ be the maximal Abelian pro-extension of~$K$. 
We assume that the  Leopoldt conjecture is verified in $K$ for all prime numbers.

\smallskip
Then there exists an Abelian extension ${\mathcal F}_\infty^1 \subseteq {\mathcal F}_\infty$ of $K$ such that ${\mathcal H}_\infty^1$ is the direct compositum over $K$ of ${\mathcal F}_\infty^1$ and the maximal $\wh \Z$-extension~$\wt {\mathcal K}_\infty$ of~$K$, and such that
$${\rm Gal}(\ov K^\ab/{\mathcal F}_\infty^1) \, \ds\mathop{\simeq}^{\rm nc} \,  \wh \Z^{r_2+1}  
\times \prd_{n \geq 1}\Big(  (\Z/2\,\Z)^\delta \times \Z/\w \, n \,\Z\Big),$$
with ${\rm Gal}(\ov K^\ab /{\mathcal H}_\infty^1) \, \ds\mathop{\simeq}^{\rm nc} \, \prd_{n \geq 1} \Big(  (\Z/2\,\Z)^\delta \times \Z/\w \, n \,\Z\Big)$. 

If $\w  = 1$, then ${\rm Gal}(\ov K^\ab/{\mathcal F}_\infty^1) \, \ds\mathop{\simeq}^{\rm nc} \,  \wh \Z^{r_2+1} \times \prd_{n \geq 1} \Z/ n \,\Z$,
 with ${\rm Gal}(\ov K^\ab /{\mathcal H}_\infty^1) \, \ds\mathop{\simeq}^{\rm nc} \, \prd_{n \geq 1} 
 \Z/ n \,\Z$.
\end{theo}

The problem for the non-imaginary quadratic fields is unchanged since in the following global exact sequence, where ${\mathcal T}_p^1 :=  {\rm Gal}(H_p^1/\wt K_p)$ for all $p$,
$$0 \to \prd_{n \geq 1} \Big( (\Z/2\,\Z)^\delta \times \Z/\w \, n \,\Z\Big) \too {\rm Gal}(\ov K^\ab/\wt {\mathcal K}_\infty)  \too \prd_p {\mathcal T}_p^1 \to 1, $$

we do not know if the structure of ${\rm Gal}(\ov K^\ab/\wt {\mathcal K}_\infty)$ can be the same for various number fields $K$ because of the unknown groups $\prd_p {\mathcal T}_p^1$ (which non-trivially depend on the $p$-adic properties of the classes and units of  the fields $K$) and the nature of the corresponding group extension.

\smallskip
 We have 
$${\rm Gal}(\ov K^\ab/K)\, \ds\mathop{\simeq}^{\rm nc} \,\wh \Z^{r_2+1} \times \prd_{n \geq 1} 
\Big( (\Z/2\,\Z)^\delta \times \Z/\w \, n \,\Z\Big)$$
as soon as the second and third condition of $p$-rationality (Remark 2.2 (i))
 are satisfied for all~$p$, which defines a weaker version of $p$-rationality which may have some interest.

\smallskip
For imaginary quadratic fields $K \ne \Q(\sqrt {-1}\,), \Q(\sqrt {-2}\,)$, we find again (since $\delta=0$
and $\w = 1$) that 
${\rm Gal}(\ov K^\ab/K)  \, \ds\mathop{\simeq}^{\rm nc} \, \wh \Z^2 \times \prd_{n \geq 1} \Z/ n \,\Z$,
as soon as, for all $p$ dividing the class number, the $p$-Hilbert class field is contained in the compositum  of the $\Z_p$-extensions\,\footnote{From [Gr1, III.2.6.6] or [Gr3, Theorem 2.3], for an imaginary quadratic field $K$,  the $2$-Hilbert class field  is contained in the compositum of its $\Z_2$-extensions if and only if $K$ is one of the 
following fields:  $\Q(\sqrt{-1}\,)$, $\Q(\sqrt{-2}\,)$, $\Q(\sqrt{-\ell}\,)$ ($\ell$ prime $\equiv 3, 5,
7 \pmod 8$), $\Q(\sqrt{-2 \ell}\,)$  ($\ell$ prime $\equiv 3, 5 \pmod8$, $\Q(\sqrt{-\ell q}\,)$ ($\ell$, $q$ primes, $\ell \equiv -q \equiv 3 \pmod 8$). For numerical studies on the groups ${\cal T}_p$, see [Cha] and [AS, \S\,7]. } 
of $K$, which is equivalent to 
$\Cl_p \, \ds\mathop{\simeq}^{\rm can} \, \Z_p \,{\rm Log}_p(I_p)\Big / \!\!\prd_{v\div p} {\rm log} (U_v^1)$.

Note that the arithmetical invariant $\prd_p {\mathcal T}_p$ (or $\prd_p {\mathcal T}_p^1$) is one of the deepest invariant of class field theory over $K$; the duality properties of each component ${\mathcal T}_p$
are related to $p$-class groups, $p$-regular kernels,\,\ldots  via reflection theorems and Galois cohomology; in the totally real case, ${\mathcal T}_p$ is connected with the $p$-adic $\zeta$-function of $K$ 
(see [Se] and [Gr1, III.2.6.5]).

\begin{rema} Let $F_p^1$ be any extension of $K$ such that $H_p^1$ is the direct compositum over $K$ of
$\wt K_p$ and $F_p^1$. From [Gr1, III.4.15.8], we know that when $F_p^1 \ne K$,
all non-trivial cyclic extensions $F_{p,i}^1 \subseteq F_p^1$ of $K$
can be embedded in a cyclic $p$-extension of arbitrarily large degree (except perhaps 
 in the special case $p=2$, $K\cap \Q(\mu^{}_{2^\infty}) = Q_\nu$, $\nu \geq 2$).

\smallskip
 Recall that the  subgroup corresponding to the compositum of the $p$-cyclically embeddable fields
(compositum which of course contains $\wt K_p$) is equal to the group $\prod_{v\div p } \mu^1_v \big /\mu_p(K)$, except perhaps in the special case where $\prod_{v\div p } \mu^1_v \big /\mu_p(K)$ may be of index $2$ in  this group. The quotient of ${\cal T}_p$ by this group is called the Bertrandias--Payan module. 

So this property shows, when $F_p^1 \ne K$, that ${\rm Gal}(\ov K^\ab{\st (p)}/H_p^1)$ cannot be a direct summand in 
${\rm Gal}(\ov K^\ab{\st (p)}/\wt K_p)$, since ${\mathcal T}_p^1$ is finite.
In other words, for any power $p^k$, taking a suitable set of cyclic extensions $F_{p,i}^1 \subseteq F_p^1$,
there exists a field $L_k \subset \ov K^\ab{\st (p)}$, such that $\wt K_p \subseteq H_p^1 \subseteq L_k$, with
${\rm Gal}(L_k/\wt K_p)$ of exponent $p^k$.
We can even assume that ${\rm Gal}(L_k/\wt K_p) \simeq (\Z/p^k\Z)^r$, where $r$ is the $p$-rank
of ${\mathcal T}_p^1$.
However, for distinct values of $k$, the fields $L_k$ may not follow any specific rule. So it is possible that only numerical computations may help to precise the structure of ${\rm Gal}(\ov K^\ab{\st (p)}/\wt K_p)$.
\end{rema}

An interesting case is that of $K= \Q(\sqrt 2\,)$ for $p=13$; in this case, $H_{13} = H_{13}^1$ is cyclic  of degree $13$ over $\wt K_{13}$ since $\varepsilon = 1 + \sqrt 2\,$ is such that $-\varepsilon^{14}$ is, modulo $13^3$, of the form $1 +13^2\,a\,\sqrt 2\,$ with $a \not \equiv 0 \pmod {13}$, which gives ${\mathcal T}_{13} \simeq \Z/13\,\Z$; indeed, use the reasonning of [Gr1, III.4.14] for real quadratic fields, or the formula given in [Gr1, III.2.6.1 (ii$_2$)] with $\Cl_{13} = 1$, $\prd_{v\div 13} U_v^1 = U_{13}^1 = 1 + 13\,(\Z_{13} \oplus \Z_{13}\,\sqrt 2\, )$.

With such similar numerical data for a real quadratic field $\Q (\sqrt d\,)$ ($p\ne 2$,  class number  prime to $p$,
$H_p^1=H_p$ of degree $p$ over $\wt K_p$, $\,\pm \varepsilon^{p+1}$ ($p$ inert) or $\pm \varepsilon^{p-1}$ ($p$ split) is, modulo $p^3$, of the form $1 + p^2\,a\,\sqrt d\,$ with a rational $a\not\equiv 0 \pmod p$), we get the following diagram:
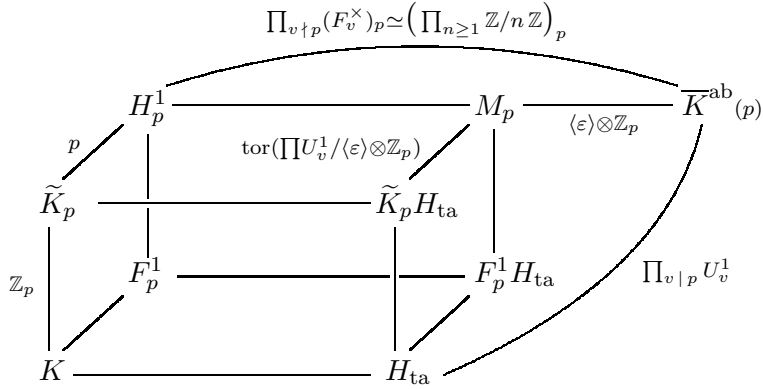
\unitlength=0.65cm
$$\vbox{\hbox{\hspace{-5.5cm}  \begin{picture}(11.5,8.0)
%
% horizontales
\put(5.5,6){\line(1,0){6}}
\put(4.0,4){\line(1,0){5.5}}

\put(5.6,2.50){\line(1,0){3}}
\put(8.6,2.50){\line(1,0){1.3}}
\put(10.1,2.50){\line(1,0){1.3}}

\put(3.5,0.50){\line(1,0){6}}

\bezier{350}(3.3,1.)(3.9,1.55)(4.5,2.1)
\bezier{350}(3.3,4.5)(3.9,5.05)(4.5,5.6)

\bezier{350}(10.3,1.)(10.9,1.55)(11.5,2.1)
\bezier{350}(10.3,4.5)(10.9,5.05)(11.5,5.6)

% vertical
\put(3.0,1.0){\line(0,1){2.5}}

\put(5.0,2.9){\line(0,1){0.98}}
\put(5.0,4.1){\line(0,1){1.3}}

\put(10.0,1.0){\line(0,1){2.5}}
\put(12.00,3.0){\line(0,1){2.5}}

\put(2.8,0.4){$K$}
\put(4.6,2.4){$F_p^1$}

\put(9.8,0.4){$H_{\rm ta}$}
\put(11.6,2.4){$F_p^1 H_{\rm ta}$}
\put(2.8,3.8){$\wt K_p$}
\put(4.6,5.8){$H_p^1$}

\put(9.6,3.8){$\wt K_p H_{\rm ta}$}
\put(11.6,5.8){$M_p$}

\put(12.6,6){\line(1,0){3}}
\put(15.8,5.8){$\ov K^\ab{\st (p)}$}

\bezier{400}(5.4,6.4)(10.45,8.0)(15.7,6.4)
\put(7.4,7.6){$\st  \prod_{ v \notdiv p} ( F_v^\times)_p \simeq \big(\prod_{n\geq 1}{\Z/n\,\Z \big)_p} $}

\put(3.4,5.1){$\st p$}

\put(2.2,2.2){$\st \Z_p$}

\bezier{300}(11,0.5)(15.5,2.5)(16.2,5.6)
\put(15.,2.5){$\st \prod_{v\div p} U_v^1$}

\put(6.9,5.){$\st {\rm tor} (\st \prod\! U_v^1/ \st  \langle\varepsilon\rangle \otimes\zz_p)$}
\put(13.5,5.5){$\st \langle\varepsilon\rangle \otimes\zz_p$}
\end{picture}   }} $$
\unitlength=1.0cm

 For $K=\Q(\sqrt 2\,)$, $p=13$, we have no more information likely to give a result on the structure
of the profinite group ${\rm Gal}(\ov K^\ab{\st (p)}/\wt K_p)$ containing a subgroup, of index $p$, 
isomorphic (non-canonically) to $\Big(\prd_{n\geq 1}{\Z/n\,\Z \Big)_p}$.

Despite the previous class field theory study, it remains possible that ${\rm Gal}(\ov K^\ab/K)$ be 
always non-canonically isomorphic to ${\wh \Z}^{r_2+1} \times \prd_{n\geq 1} 
\Big( (\Z/2\,\Z)^\delta \times \Z/\w \, n \,\Z\Big)$, independently of  additional arithmetic 
considerations about the group $\prod_{p}{\mathcal T}_p^1$.
If not (more probable), a description of the profinite group ${\rm Gal}(\ov K^\ab/K)$ may be very tricky.
Any information will be welcome.

\bigskip
\centerline{\bf Acknowledgements}

\medskip
I thank Peter Stevenhagen who pointed out to me an error regarding the definition of the parameter $\w$ associated to the number field $K$ in my previous versions (arXiv:1212.3588). This error also exists in [O] and in the first draft (arXiv:1209.6005) of [AS] but will be corrected in the final paper in the proceedings volume of ANTS-X, San Diego 2012. I also thank him for his remarks and cooperation in the improvement of this version.

\end{document}